\documentclass[12pt]{article}
\usepackage {amsmath, amsthm, amssymb, mathrsfs}
\usepackage{enumerate}
\usepackage[usenames]{color}
\usepackage{indentfirst}
\newtheorem{thm}{Theorem}
\newtheorem{lem}[thm]{Lemma}

\newtheorem*{C}{Claim}
\newtheorem{e}[thm]{Example}
\newcommand{\ex}[1]{\begin{e}\upshape #1 \end{e}}
\newtheorem*{Q}{Question}
\newtheorem*{p}{Proof}
\newcommand{\pf}[1]{\begin{p}\upshape #1 \qed \end{p}}
\newtheorem*{pC}{Proof of Claim}
\newcommand{\pfC}[1]{\begin{pC}\upshape #1 \qed \end{pC}}
\newtheorem*{re}{Remark}
\newcommand{\rem}[1]{\begin{re}\upshape #1 \end{re}}
\let\svthefootnote\thefootnote
\title{{Compact leaves of the foliation defined by the kernel of a $T^2$-invariant presymplectic form}
\let\thefootnote\relax\footnotetext{\hspace{-20pt}
2020 {\it Mathematics Subject Classification.} Primary 37C86, Secondary 37C85, 53D10\\
{\it Key Words and Phrases.} presymplectic form, foliation, torus action, generalized Weinstein conjecture
}\let\thefootnote\svthefootnote}
\author{{Asuka Hagiwara}
\footnote{Graduate School of Science and Engineering, Ibaraki University, Bunkyo 2-1-1, Mito, Ibaraki, 310-8512, Japan.
}}
\newcommand{\Addresses}{{
  \bigskip
\noindent
{\it Present Address}:\\
Numata Girls' High School,\\
Higashikurauchimachi 753-3, Numata, Gumma, 378-0043, Japan.\\
{\it e-mail}: hagiwara-asuka@edu-g.gsn.ed.jp
}}
\date{ }
\begin{document}
\maketitle

\begin{abstract}
We investigate the foliation defined by the kernel of an exact presymplectic form $d\alpha$ of rank $2n$ on a $(2n+r)$-dimensional closed manifold $M$.
For $r=2$, we prove that the foliation has at least two leaves which are homeomorphic to a $2$-dimensional torus, if $M$ admits a locally free $T^2$-action which preserves $d\alpha$ and satisfies that the function $\alpha(Z_2)$ is constant, where $Z_1,Z_2$ are the infinitesimal generators of the $T^2$-action.
We also give its generalization for $r\geq1$.
\end{abstract}

\section{Introduction}

Let $n$ and $r$ be positive integers.
A smooth manifold $M$ of dimension $2n+r$ is said to be {\it presymplectic} if it carries a closed $2$-form $\omega$ such that the dimension of the kernel of $\omega(x)$ is $r$ for all $x\in M$, that is, $\omega$ is of constant rank $2n$.
The $2$-form $\omega$ is called a {\it presymplectic form} on $M$.
We denote the kernel of $\omega$ by $\mbox{ker}\, \omega$, which is an involutive distribution of dimension $r$.
Hence, by Frobenius' theorem, $\mbox{ker}\, \omega$ defines an $r$-dimensional foliation $\mathcal{F}_{\omega}$.
The leaves of $\mathcal{F}_{\omega}$ are integral manifolds of the distribution. 

A special case of presymplectic manifolds is when a manifold $M$ is $(2n+1)$-dimensional and carries a contact form $\alpha$, that is, $\alpha$ is a $1$-form on $M$ such that $\alpha \wedge (d\alpha)^n\neq0$.
In this case, the $2$-form $\omega:=d\alpha$ is a presymplectic form on $M$ and the foliation $\mathcal{F}_{\omega}$ is $1$-dimensional.
The leaves of $\mathcal{F}_{\omega}$ are integral curves of $\mbox{ker}\, \omega$, which are called {\it characteristics}.
The Weinstein conjecture \cite{Weinstein1979} asserts that there exists at least one closed characteristic on a closed contact manifold, which is one of the fundamental problem concerning dynamics on contact manifolds.
This conjecture is still open, but has been proved in several cases.
For instance, Viterbo  \cite{Viterbo1987} proved it for compact contact manifolds which are hypersurfaces of contact type in $(\mathbb{R}^{2n},\omega_0)$, where $\omega_0=\sum_{j=1}^{n} dy_j \wedge dx_j$ is the standard symplectic form.
Hofer \cite{Hofer1993} proved this conjecture for $S^3$ and for overtwisted $3$-dimensional contact manifolds, and after that Taubes \cite{Taubes2007} solved it affirmatively in dimension $3$.

In \cite{B-R. 1995}, Banyaga and Rukimbira raised the following question for presymplectic manifolds $(M,\omega)$, which generalizes the Weinstein conjecture.

\begin{Q}[{\cite[p.$3902$]{B-R. 1995}}]\upshape
When does the foliation $\mathcal{F}_{\omega}$ admits a compact leaf?
\end{Q}

In addition, they focused their attention on the case where $\mbox{dim}\, \mathcal{F}_{\omega}=1$ and proved the following theorem.

\begin{thm}[Banyaga-Rukimbira, \cite{B-R. 1995}]\label{T.B-R}
Let $M$ be a $(2n+1)$-dimensional oriented closed $C^{\infty}$-manifold with a $1$-form $\alpha$ such that the $2$-form $\omega:=d\alpha$ has constant rank $2n$ everywhere.
If there exists a locally free circle-action on $M$ which preserves $\omega$, then $\mathcal{F}_{\omega}$ has at least two closed leaves.
\end{thm}

In particular, they pointed out that this result guarantees the existence of periodic orbits of the Hamiltonian system on hypersurfaces in $\mathbb{R}^{2n}$ which are not necessarily of contact type.
However, they did not mention the case where $\mbox{dim}\, \mathcal{F}_{\omega}\geq2$.
For such a direction, it seems that there are few works which assume that $(M,\omega)$ is $r$-contact
(see Remark in Section $3$, and \cite{Finamore} for details).

In this paper, we consider a higher dimensional version of Theorem \ref{T.B-R} without assuming that $(M,\omega)$ is $r$-contact.
We show the following result for the case where $\mbox{dim}\, \mathcal{F}_{\omega}=2$.

\begin{thm}\label{T.main}
Let $M$ be a $(2n+2)$-dimensional oriented closed $C^{\infty}$-manifold with an
exact presymplectic form $\omega$, that is, there exists a $1$-form $\alpha \in \Omega^1(M)$ such that the $2$-form $\omega=d\alpha$ has rank $2n$ everywhere.
Assume that $M$ admits a locally free $T^2$-action with the following conditions:
\vspace{-10pt}
\begin{enumerate}[{\rm (i)}]
\setlength{\itemsep}{-5pt}
\item The $T^2$-action preserves $\omega$,
\item The function $\alpha(Z_2)$ is constant on $M$,
\end{enumerate}\vspace{-10pt}
where $Z_1,Z_2$ denote the infinitesimal generators of the $T^2$-action.
Then
\vspace{-10pt}
\begin{enumerate}[{\rm (1)}]
\setlength{\itemsep}{-5pt}
\setlength{\leftskip}{6.5pt}
\item the $2$-dimensional $C^{\infty}$-foliation $\mathcal{F}_{\omega}$ has at least two leaves which are homeomorphic to a $2$-dimensional torus.
\item Moreover, if the function $\alpha(Z_1)$ is also constant, then $\mathcal{F}_{\omega}$ coincides with the foliation of the $T^2$-action and hence, all leaves of $\mathcal{F}_{\omega}$ are homeomorphic to a $2$-dimensional torus.
\end{enumerate}\vspace{-10pt}
\end{thm}

We prove Theorem \ref{T.main} in the next section.
In Section $3$ we provide examples of Theorem \ref{T.main} and mention that a similar result holds for the general case where $\mbox{dim}\, \mathcal{F}_{\omega}=r$.

\section{Proof of Theorem 2}
%
Let $\rho:M\times T^2 \to M$ be a locally free $T^2$-action on $M$ with the conditions (i) and (ii) above, and let $\mathcal{D}:M\ni q\mapsto \mathcal{D}_q=\mbox{span}\{ Z_1(q), Z_2(q) \} \subset T_qM$ be the distribution determined by $\rho$, where $Z_1,Z_2$ denote the infinitesimal generators of $\rho$.
The distribution $\mathcal{D}$ is involutive and of dimension $2$ since $\rho$ is locally free. Therefore, now $M$ has two foliations of codimension $2n$ which are determined by $\mbox{ker}\, \omega$ and $\mathcal{D}$.
Note that, by definition, every leaf of $\mathcal{D}$ is homeomorphic to a $2$-dimensional torus $T^2$.

We define the diffeomorphism $s:M\to M$ by $s(x):=\rho(x,s)$ for $s\in T^2$ and the $1$-form
$$
\alpha_0:=\int_{T^2} (s^*\alpha)\, d\sigma,
$$
where $\sigma$ is the Haar measure on $T^2$.
Then $\alpha_0$ is $T^2$-invariant.
Thus,
$$
L_{Z_i}\alpha_0=\lim_{t\to 0}\frac{(\varphi_i^t)^*\alpha_0-\alpha_0}{t}=0,
\quad i=1,2,
$$
where $\varphi_i^t=\rho(\cdot,\exp tZ_i)$ denotes the flow of $Z_i$, respectively.
Moreover, define smooth functions $S_i:M\to \mathbb{R}$ by
\begin{equation*}\label{S_1}
S_i(x) :=-\alpha_0(x)(Z_i(x))= -i_{Z_i}\alpha_0(x),
\quad i=1,2.
\end{equation*}
By the condition (i), we have
$$
d\alpha_0=\int_{T^2} d(s^*\alpha) \, d\sigma=\int_{T^2} (s^*d\alpha) \, d\sigma=\omega,
$$
and hence, due to the Cartan formula,
\begin{equation}\label{dS_i}
dS_i=-di_{Z_i}\alpha_0=-L_{Z_i}\alpha_0+i_{Z_i}d\alpha_0=i_{Z_i}\omega.
\end{equation}
Since the group $T^2$ is commutative, we have $s\, \circ \, \varphi_2^t(x)=\varphi_2^t \, \circ \, s(x)$ for all $x\in M$.
Differentiating this equation in $t$, we obtain $ds(x)Z_2(x)=Z_2(s(x)),$ so that
$$
S_2(x)= -\int_{T^2} s^*\alpha(x)(Z_2(x)) \, d\sigma= -\int_{T^2} \alpha(s(x))(Z_2(s(x))) \, d\sigma.
$$
Thus, by the condition (ii) the function $S_2$ is constant and so $dS_2=0$.
Hence,
\begin{equation}\label{Z_2inker}
Z_2(q)\in \mbox{ker}\, \omega(q),\quad \forall q\in M 
\end{equation}
by (\ref{dS_i}) for $i=2$.
To investigate the relation between $Z_1$ and $\mbox{ker}\, \omega$, we shall take a special chart.
Since $\mbox{ker}\, \omega$ is involutive, by using Frobenius' theorem \cite[p.89, Theorem 1]{Chevalley}, for any point $p\in M$ we can choose the following chart $(\varphi,U)$ around $p$:
for $q\in U$ we write
$$
\varphi(q)=(x(q),z_1(q),z_2(q)),\quad x(q)=(x_1(q),\ldots,x_{2n}(q)),
$$
then this chart satisfies that $\varphi(p)=(0,0,0)$ and
$$
\xi_i:=\frac{\partial}{\partial z_i},
\quad i=1,2
$$
form a local frame of $\mbox{ker}\, \omega$ on $U$.
For $q\in U$ we put
$$
Z_1(q)
=\sum_{j=1}^{2n} X_j(x,z_1,z_2)\frac{\partial}{\partial x_j}
+u_1(x,z_1,z_2)\frac{\partial}{\partial z_1}
+u_2(x,z_1,z_2)\frac{\partial}{\partial z_2}
$$
and prove the following

\begin{lem}\label{X_j(x)}\label{L.1}
The function $X_j$ is independent of $z_1,z_2$ for any $j=1,\ldots,2n$.
\end{lem}

\pf{
In $U$, since $\xi_i \in \mbox{ker}\, \omega$, we have
$
L_{\xi_i}\omega = i_{\xi_i}d\omega + di_{\xi_i}\omega = 0
$
and
\begin{align*}
L_{\xi_i}(i_{Z_j}\omega)
&=i_{\xi_i}d(i_{Z_j}\omega) + di_{\xi_i}(i_{Z_j}\omega)\\
&=i_{\xi_i}ddS_j - di_{Z_j}i_{\xi_i}\omega= 0
\quad(\because \ (\ref{dS_i})).
\end{align*}
Therefore,
\begin{equation*}\label{i[xi,Z]}
i_{[\xi_i, Z_j]}\omega = L_{\xi_i}i_{Z_j}\omega -i_{Z_j}L_{\xi_i}\omega = 0
\end{equation*}
and namely, there exist functions $\lambda_i^1, \lambda_i^2:U \to \mathbb{R}$ such that
\begin{equation*}\label{[Z_i]}
\left[ \xi_i,Z_1 \right]=\lambda_i^1 \xi_1 + \lambda_i^2 \xi_2.
\end{equation*}
On the other hand,
\begin{align*}
[\xi_i,Z_1]
=
\left[ \frac{\partial}{\partial z_i}, Z_1 \right]
&=
\left[ \frac{\partial}{\partial z_i},\sum_{j=1}^{2n} X_j \frac{\partial}{\partial x_j} \right]
+\left[ \frac{\partial}{\partial z_i},u_1 \frac{\partial}{\partial z_1} \right]
+\left[ \frac{\partial}{\partial z_i},u_2 \frac{\partial}{\partial z_2} \right]\\
&=
\sum_{j=1}^{2n} \frac{\partial X_j}{\partial z_i} \frac{\partial}{\partial x_j}
+\frac{\partial u_1}{\partial z_i} \frac{\partial}{\partial z_1}
+\frac{\partial u_2}{\partial z_i} \frac{\partial}{\partial z_2}.
\end{align*}
Therefore,
$
\partial X_j/\partial z_i = 0.
$}

By this lemma, we have
$$
Z_1(q)
=\sum_{j=1}^{2n} X_j(x)\frac{\partial}{\partial x_j}
+u_1(x,z_1,z_2)\frac{\partial}{\partial z_1}
+u_2(x,z_1,z_2)\frac{\partial}{\partial z_2}
\quad
\mbox{for}\ q\in U.
$$

Now, we assume that $p\in M$ is a critical point of $S_1$.
By (\ref{dS_i}) for $i=1$, we have $Z_1(p)\in \mbox{ker}\, \omega(p)$.
Because $\varphi(p) = (0,0,0)$, we have
$$
Z_1(p)
=\sum_{j=1}^{2n} X_j(0)\frac{\partial}{\partial x_j}
+u_1(0,0,0)\frac{\partial}{\partial z_1}
+ u_2(0,0,0)\frac{\partial}{\partial z_2},
$$
so that
\begin{equation}\label{X(0)}
X_j(0) =0,\quad j=1,\ldots ,2n.
\end{equation}
We denote by $\mathcal{F}_{\omega}(p)$ the leaf of $\mathcal{F}_{\omega}$ passing through the point $p$.
Let $q\in \mathcal{F}_{\omega}(p) \cap U$.  
Then there exist $t_1,t_2\in \mathbb{R}$ such that $\varphi(q) = (0,t_1,t_2)$, since
$$
T_q(\mathcal{F}_{\omega}(p)) = \mbox{ker}\, \omega(q) = \mbox{span}\left\{ \left( \frac{\partial}{\partial z_1} \right)_q,\left( \frac{\partial}{\partial z_2} \right)_q \right\}.
$$
Thus, by (\ref{X(0)}), we obtain
\begin{equation}\label{Z_1inker}
\begin{aligned}
Z_1(q)
&=\sum_{j=1}^{2n} X_j(0)\frac{\partial}{\partial x_j}
+u_1(0,t_1,t_2)\frac{\partial}{\partial z_1}
+ u_2(0,t_1,t_2)\frac{\partial}{\partial z_2}\\
&=u_1(0,t_1,t_2)\frac{\partial}{\partial z_1}
+ u_2(0,t_1,t_2)\frac{\partial}{\partial z_2}\in \mbox{ker}\, \omega(q).
\end{aligned}
\end{equation}
It follows from (\ref{Z_2inker}) that
\begin{equation}\label{ker=D1}
\mbox{ker}\, \omega(q)=\mathcal{D}_q,\quad \forall q\in \mathcal{F}_{\omega}(p)\cap U,
\end{equation}
since the dimensions of $\mbox{ker}\, \omega$ and $\mathcal{D}$ are the same.
In order to verify that $\mbox{ker}\, \omega$ coincides with $\mathcal{D}$ on the leaf $\mathcal{F}_{\omega}(p)$, we cover $\mathcal{F}_{\omega}(p)$ by a set of charts $\{(U_a,\varphi_a)\}_{a\in A}$, where each $(U_a,\varphi_a)$ is a chart as we chose above.
For $(p\neq)\, \hat{p}\in \mathcal{F}_{\omega}(p)$ we take a chart $(\hat{U},\hat{\varphi})\in \{(U_a,\varphi_a)\}_{a\in A}$ around $\hat{p}$ which satisfies $\mathcal{F}_{\omega}(p)\cap U\cap \hat{U}\neq \emptyset$. 
We write
$$
\hat{\varphi}(q)=(y(q),w_1(q),w_2(q)),\quad y(q)=(y_1(q),\ldots,y_{2n}(q))
$$
for $q\in \hat{U}$, then $\hat{\varphi}(\hat{p})=(0,0,0)$ and
$
\partial/\partial w_1,\partial/\partial w_2
$
form a local frame of $\mbox{ker}\, \omega$ on $\hat{U}$.
With this setting, by applying Lemma \ref{L.1}, we have
\begin{equation*}
Z_1(q)
=\sum_{j=1}^{2n} \hat{X}_j(y)\frac{\partial}{\partial y_j}
+\hat{u}_1(y,w_1,w_2)\frac{\partial}{\partial w_1}
+ \hat{u}_2(y,w_1,w_2)\frac{\partial}{\partial w_2}
\end{equation*}
in $\hat{U}$ as well.
Since $T_{q}(\mathcal{F}_{\omega}(p))=\mbox{ker}\, \omega(q)$, we have $\hat{\varphi}(q)=(0,w_1(q),w_2(q))$ for all $q\in \mathcal{F}_{\omega}(p)\cap \hat{U}$.
Therefore, for $q\in \mathcal{F}_{\omega}(p)\cap U\cap \hat{U}$,
$$
Z_1(q)
=\sum_{j=1}^{2n} \hat{X}_j(0)\frac{\partial}{\partial y_j}
+\hat{u}_1(0,w_1(q),w_2(q))\frac{\partial}{\partial w_1}
+ \hat{u}_2(0,w_1(q),w_2(q))\frac{\partial}{\partial w_2}
$$
holds.
From (\ref{ker=D1}), we obtain $\hat{X}_j(0)=0$, $j=1,\ldots,2n$.
Thus, for all $q\in \mathcal{F}_{\omega}(p)\cap \hat{U}$ we have
\begin{equation*}
Z_1(q)
=\hat{u}_1(0,w_1(q),w_2(q))\frac{\partial}{\partial w_1}
+ \hat{u}_2(0,w_1(q),w_2(q))\frac{\partial}{\partial w_2}\in \mbox{ker}\, \omega(q),
\end{equation*}
so that $\mbox{ker}\, \omega(q)=\mathcal{D}_{q}$.
Repeating this argument we see that
\begin{equation}\label{ker=D}
\mbox{ker}\, \omega(q)=\mathcal{D}_q,\quad \forall q\in \mathcal{F}_{\omega}(p).
\end{equation}
Therefore, $\mathcal{F}_{\omega}(p)$ is not only an integral manifold of $\mbox{ker}\, \omega$, but also that of $\mathcal{D}$ passing through $p$.
On the other hand, the $T^2$-orbit $T^2(p)$ of $p$ is nothing but the maximal connected integral manifold of $\mathcal{D}$ passing through $p$.
Hence, by \cite[p.172, Theorem 1]{Matsushima}, $\mathcal{F}_{\omega}(p)$ is an open submanifold in $T^2(p)$.
Actually, we claim the following

\begin{lem}\label{F=T}
$\mathcal{F}_{\omega}(p)=T^2(p)$ holds.
In particular, $\mathcal{F}_{\omega}(p)$ is homeomorphic to $T^2$.
\end{lem}

\pf{
Since $\mathcal{F}_{\omega}(p)$ is a nonempty open set in a connected space $T^2(p)$, it suffices to show that $\mathcal{F}_{\omega}(p)$ is a closed set in $T^2(p)$.
We denote by $\overline{\mathcal{F}_{\omega}(p)}$ the closure of $\mathcal{F}_{\omega}(p)$ in $T^2(p)$.
Let $q_0\in \overline{\mathcal{F}_{\omega}(p)}$.
We take the following chart $(U,\varphi)$ around $q_0$ as we chose above:
for $q\in U$ we write
$$
\varphi(q)=(x(q),z_1(q),z_2(q)),\quad x(q)=(x_1(q),\ldots,x_{2n}(q)),
$$
then $\varphi(q_0)=(0,0,0)$, and
$
\partial/\partial z_1, \partial/\partial z_2
$
form a local frame of $\mbox{ker}\, \omega$ on $U$.
By Frobenius' theorem \cite[p.89, Theorem 1]{Chevalley}, the slice
$$
S_0:=\{q\in U\mid x_1(q)=0,\ldots,x_{2n}(q)=0\}
$$
is a connected integral manifold of $\mbox{ker}\, \omega$ passing through $q_0$.

\begin{C}\label{C.1}\upshape
The slice $S_0$ is an open submanifold in $T^2(p)$ which contains $q_0$.
\end{C}

\pfC{
For $q\in U$, by Lemma \ref{L.1}, we have
\begin{equation*}\label{Z_1}
Z_1(q)
=\sum_{j=1}^{2n} X_j(x)\frac{\partial}{\partial x_j}
+u_1(x,z_1,z_2)\frac{\partial}{\partial z_1}
+ u_2(x,z_1,z_2)\frac{\partial}{\partial z_2}.
\end{equation*}
First we prove that $Z_1(q_0)\in \mbox{ker}\, \omega(q_0)$.
Arguing by contradiction we assume that $Z_1(q_0)\notin \mbox{ker}\, \omega(q_0)$, that is, $X_j(x(q_0))=X_j(0)\neq0$ for some $j\in \{ 1,\ldots,2n\}$.
Since the function $X_j$ is smooth on $U$, there exists a neighborhood $U'\subset U$ of $q_0$ such that
$
X_j(x(q))\neq0
$
for any $q\in U'$.
Because $q_0\in \overline{\mathcal{F}_{\omega}(p)}$, the intersection $\mathcal{F}_{\omega}(p)\cap U'$ is nonempty, so that for $q\in \mathcal{F}_{\omega}(p)\cap U'$ it holds that
$
X_j(x(q))\neq0.
$
This is a contradiction to (\ref{ker=D}).
Hence, we have $X_j(0)=0$ for all $j=1,\ldots,2n$.

Therefore, for all $q\in S_0$, we can write $\varphi(q)=(0,z_1(q),z_2(q))$ and we have
\begin{equation*}
Z_1(q)
=u_1(0,z_1(q),z_2(q))\frac{\partial}{\partial z_1}
+ u_2(0,z_1(q),z_2(q))\frac{\partial}{\partial z_2}\in \mbox{ker}\, \omega(q).
\end{equation*}
Thus, by (\ref{Z_2inker}) we obtain $\mathcal{D}_q=\mbox{ker}\, \omega(q)$ for all $q\in S_0$.
It follows from $T_q S_0=\mbox{ker}\, \omega(q)$ that $S_0$ is a connected integral manifold of $\mathcal{D}$ through $q_0$.
On the other hand, $T^2(p)$ is the maximal connected integral manifold of $\mathcal{D}$ which contains $q_0$.
Therefore, by \cite[p.172, Theorem 1]{Matsushima}, we see that $S_0$ is an open submanifold in $T^2(p)$ which contains $q_0$.
}
By this claim, $S_0$ is an open neighborhood of $q_0$ in $T^2(p)$.
Because $q_0\in \overline{\mathcal{F}_{\omega}(p)}$, the intersection $S_0\cap \mathcal{F}_{\omega}(p)$ is nonempty.
Thus, $\mathcal{F}_{\omega}(p)$ is the leaf of $\mbox{ker}\, \omega$ passing through a point $p_0\in S_0\cap \mathcal{F}_{\omega}(p)$ as well.
Due to Frobenius' theorem \cite[p.89, Theorem 1]{Chevalley}, $S_0$ is a connected integral manifold of $\mbox{ker}\, \omega$ which contains the point $p_0$.
Hence, by \cite[p.172, Theorem 1]{Matsushima}, $S_0$ is also an open submanifold in $\mathcal{F}_{\omega}(p)$.
Consequently, $q_0\in S_0 \subset \mathcal{F}_{\omega}(p)$ and so we have $\overline{\mathcal{F}_{\omega}(p)}=\mathcal{F}_{\omega}(p)$.
}

Finally, we check the existence of critical points of the function $S_1$ on $M$.
If $S_1$ is not constant, then it has at least two critical points $p_{\rm{max}},p_{\rm{min}}\in M$ corresponding to a maximum and a minimum, respectively, since $M$ is closed.
Moreover, $p_{\rm{max}}$ and $p_{\rm{min}}$ are on different leaves of $\mathcal{F}_{\omega}$.
This follows from the following.
By (\ref{dS_i}) and (\ref{Z_2inker}), we have
$$
L_{Z_1}S_1 = i_{Z_1}dS_1 + di_{Z_1}S_1 = i_{Z_1}i_{Z_1}\omega = 0
$$
and
$$
L_{Z_2}S_1 = i_{Z_2}dS_1 + di_{Z_2}S_1 = i_{Z_2}i_{Z_1}\omega = -i_{Z_1}i_{Z_2}\omega = 0,
$$
so that $S_1$ is constant along $T^2(p)=\{ \varphi_1^s \, \circ \, \varphi_2^t(p) \mid s,t\in \mathbb{R} \}$.
Therefore, $T^2(p_{\rm{max}}) \cap T^2(p_{\rm{min}}) = \emptyset$ and hence, by Lemma \ref{F=T}, we obtain $\mathcal{F}_{\omega}(p_{\rm{max}})\cap \mathcal{F}_{\omega}(p_{\rm{min}})=\emptyset$.
If $S_1$ is constant on $M$, then $dS_1(p)=0$ for all $p\in M$, so that $\mbox{ker}\, \omega=\mathcal{D}$ on $M$.
Hence, from the uniqueness of foliations, $\mathcal{F}_{\omega}$ coincides with the foliation defined by $\mathcal{D}$ and therefore, each leaf of $\mathcal{F}_{\omega}$ is homeomorphic to $T^2$.
Consequently, in any case, we see that $\mathcal{F}_{\omega}$ has at least two leaves which are homeomorphic to $T^2$.
In particular, if $\alpha(Z_1)$ is constant, then, by a similar calculation for $S_2$, we see that $S_1$ is constant on $M$.
Thus, we complete the proof of Theorem \ref{T.main}.

\section{Examples and Remarks}
%
We first provide an example of Theorem 2.

\ex{\label{E.1}
We consider the case of a submanifold of codimension 2 of a symplectic manifold $(\mathbb{R}^6,d\lambda)$, where
$$
\lambda=\frac{1}{2}\sum_{j=1}^3 (y_j dx_j -x_j dy_j)
$$
is the Liouville form.
Then $d\lambda=\omega_0$ is the standard symplectic form on $\mathbb{R}^6$.
We define functions $G_i:\mathbb{R}^6\to \mathbb{R}$, $i=1,2$, by
\begin{equation*}
G_1:=(x_1^2+y_1^2)^2+x_2^2+y_2^2,
\quad
G_2:=x_2^2+y_2^2 + x_3^2+y_3^2
\end{equation*}
and for real numbers $c_1>c_2>0$ we put
\begin{equation*}
M:=\{(x,y)\in \mathbb{R}^6\mid G_1=c_1,\ G_2=c_2\}.
\end{equation*}
Since the gradients
\begin{equation*}
\nabla G_1=
\begin{pmatrix}
4x_1(x_1^2+y_1^2)\\
4y_1(x_1^2+y_1^2)\\
2x_2\\
2y_2\\
0\\
0
\end{pmatrix},\quad
\nabla G_2=
\begin{pmatrix}
0\\
0\\
2x_2\\
2y_2\\
2x_3\\
2y_3
\end{pmatrix}\smallskip
\end{equation*}
are linearly independent on $M$, we see that $M$ is a $4$-dimensional closed submanifold of $\mathbb{R}^6$.
Because $\{G_1,G_2\}=0$, we deduce that $\mbox{dim}(\mbox{ker}\, (\omega_0|_M))=2$ and $\omega_0|_M$ is a presymplectic form of constant rank $2$ (see \cite[p.27]{Moser1978}).
We define a $T^2$-action $\rho$ on $M$ as follows.
We put $z_j=x_j+\sqrt{-1}\, y_j$ and for $(e^{\sqrt{-1}\,2\pi s_1},e^{\sqrt{-1}\,2\pi s_2})\in T^2$, $s_1,s_2\in \mathbb{R}$, we set
\begin{equation*}
f_{s_1}:
\begin{pmatrix}
z_1\\
z_2\\
z_3
\end{pmatrix}
\longmapsto
\begin{pmatrix}
e^{\sqrt{-1}\,2\pi s_1}&0&0\\
0&e^{\sqrt{-1}\,2\pi s_1}&0\\
0&0&e^{\sqrt{-1}\,2\pi s_1}
\end{pmatrix}
\begin{pmatrix}
z_1\\
z_2\\
z_3
\end{pmatrix},
\end{equation*}
\begin{equation*}
g_{s_2}:
\begin{pmatrix}
z_1\\
z_2\\
z_3
\end{pmatrix}
\longmapsto
\begin{pmatrix}
1&0&0\\
0&e^{\sqrt{-1}\,2\pi s_2}&0\\
0&0&e^{\sqrt{-1}\,2\pi s_2}
\end{pmatrix}
\begin{pmatrix}
z_1\\
z_2\\
z_3
\end{pmatrix}
\end{equation*}
and define
$$
\rho^{(e^{\sqrt{-1}\,2\pi s_1},e^{\sqrt{-1}\,2\pi s_2})}(z):=f_{s_1} \circ g_{s_2} (z)=g_{s_2} \circ f_{s_1} (z),
\quad z=(z_1,z_2,z_3)\in M.
$$
Then $\rho^{(e^{\sqrt{-1}\,2\pi s_1},e^{\sqrt{-1}\,2\pi s_2})}(z)\in M$ and we can easily check that
\begin{equation}\label{lambda}
{\rho^{(e^{\sqrt{-1}\,2\pi s_1},e^{\sqrt{-1}\,2\pi s_2})}}^* \lambda=\lambda
\end{equation}
and $\rho$ is a free action preserving $d\lambda$.
The infinitesimal generators of $\rho$ are given by
\begin{align*}
Z_1
&=
\left. \frac{d}{ds_1}f_{s_1}(x,y)\right|_{s_1=0}
=
\sum_{j=1}^3 \left(-2\pi y_j\frac{\partial}{\partial x_j}+2\pi x_j\frac{\partial}{\partial y_j}\right),\\
Z_2
&=
\left. \frac{d}{ds_2} g_{s_2}(x,y)\right|_{s_2=0}
=
\sum_{j=2}^3 \left(-2\pi y_j\frac{\partial}{\partial x_j}+2\pi x_j\frac{\partial}{\partial y_j}\right),
\end{align*}
and therefore, the function
\begin{equation*}
\lambda(Z_1)
=-\pi \sum_{j=1}^3(x_j^2+y_j^2)
=-\pi (x_1^2+y_1^2+c_2)
\end{equation*}
is nonconstant on $M$ and the function
\begin{equation*}
\lambda(Z_2)=-\pi (x_2^2+y_2^2+x_3^2+y_3^2)=-\pi c_2
\end{equation*}
is constant on $M$.
It follows that $(M,\omega_0|_M)$ satisfies the conditions (i), (ii) of Theorem \ref{T.main}, but not the additional condition in (2) of Theorem \ref{T.main}.
Thus, the foliation $\mathcal{F}_{\omega_0|_M}$ has at least two leaves which are homeomorphic to $T^2$.

In this case, by (\ref{lambda}), the function $S_1$ defined in the proof of Theorem \ref{T.main} is given by
$$
S_1=-\int_{T^2} \left({\rho^{(e^{\sqrt{-1}\,2\pi s_1},e^{\sqrt{-1}\,2\pi s_2})}}^* \lambda \right)(Z_1)\, d\sigma
=-\int_{T^2} \lambda(Z_1) \, d\sigma 
=-\lambda(Z_1).
$$
Due to the compactness of $M$, there exist $p_{\rm{max}},p_{\rm{min}}\in M$ such that $S_1(p_{\rm{max}})$ is a maximal value and $S_1(p_{\rm{min}})$ is a minimal value of $S_1$.
By the former part of the last paragraph of the proof of Theorem \ref{T.main}, the function $S_1$ is constant along the leaves $\mathcal{F}_{\omega_0|_M}(p_{\rm{max}})$ and $\mathcal{F}_{\omega_0|_M}(p_{\rm{min}})$.
In fact, by using the method of Lagrange multipliers, we see thet $S_1$ has only two critical submanifolds and $S_1$ takes the maximal value $S_1(p_{\rm{max}})=\pi (\sqrt{c_1}+c_2)$ on
$$
\mathcal{F}_{\omega_0|_M}(p_{\rm{max}})
=\{ (x,y)\in M\mid x_1^2+y_1^2=\sqrt{c_1},\ x_2^2=y_2^2=0,\ x_3^2+y_3^2=c_2\}
$$
and the minimal value $S_1(p_{\rm{min}})=\pi (\sqrt{c_1-c_2}+c_2)$ on
$$
\mathcal{F}_{\omega_0|_M}(p_{\rm{min}})
=\{ (x,y)\in M\mid x_1^2+y_1^2=\sqrt{c_1-c_2},\ x_2^2+y_2^2=c_2,\ x_3=y_3=0\}.
$$
Moreover, we can easily check that
\begin{equation}\label{eq.Z_1}
Z_1(q)
\notin \mbox{ker}\, \omega(q)
\quad
\mbox{for}\ q\in M\setminus (\mathcal{F}_{\omega_0|_M}(p_{\rm{max}}) \cup \mathcal{F}_{\omega_0|_M}(p_{\rm{min}})).
\end{equation}}

Another example can be found in \cite[p.24]{Moser1978}, which satisfies the condition in $(2)$
of Theorem 2.

\ex{\label{E.2}
We also consider $(\mathbb{R}^6,\omega_0=d\lambda)$, where $\lambda$ is the Liouville form.
Define functions on $\mathbb{R}^6$ by
\begin{equation*}
G_1:=x_1^2+y_1^2,
\quad
G_2:=\sum_{j=1}^3 (x_j^2+y_j^2).
\end{equation*}
For real numbers $c_2>c_1>0$, we set
\begin{equation*}
M:=\{(x,y)\in \mathbb{R}^6\mid G_1=c_1,\ G_2=c_2\}
\end{equation*}
then, as in Example \ref{E.1}, $M$ is a $4$-dimensional closed manifold with a presymplectic form $\omega_0|_M$ of constant rank $2$.
We shall define a $T^2$-action on $M$.
Put $z_j=x_j+\sqrt{-1}\, y_j$.
For $(e^{\sqrt{-1}\,2\pi s_1},e^{\sqrt{-1}\,2\pi s_2})\in T^2$, $s_1,s_2\in \mathbb{R}$, we set
\begin{equation*}
f_{s_1}:
\begin{pmatrix}
z_1\\
z_2\\
z_3
\end{pmatrix}
\longmapsto
\begin{pmatrix}
e^{\sqrt{-1}\,2\pi s_1}&0&0\\
0&1&0\\
0&0&1
\end{pmatrix}
\begin{pmatrix}
z_1\\
z_2\\
z_3
\end{pmatrix}
,
\end{equation*}
\begin{equation*}
g_{s_2}:
\begin{pmatrix}
z_1\\
z_2\\
z_3
\end{pmatrix}
\longmapsto
\begin{pmatrix}
e^{\sqrt{-1}\,2\pi s_2}&0&0\\
0&e^{\sqrt{-1}\,2\pi s_2}&0\\
0&0&e^{\sqrt{-1}\,2\pi s_2}
\end{pmatrix}
\begin{pmatrix}
z_1\\
z_2\\
z_3
\end{pmatrix}
\end{equation*}
and define
$$
\rho^{(e^{\sqrt{-1}\,2\pi s_1},e^{\sqrt{-1}\,2\pi s_2})}(z)
:=f_{s_1} \circ g_{s_2}(z)
=g_{s_2} \circ f_{s_1}(z),
\quad z=(z_1,z_2,z_3)\in M.
$$
We see that this action is free and preserves $\omega_0$.
Similar to Example \ref{E.1}, we denote the infinitesimal generators of $\rho$ by $Z_1,Z_2$.
Then the functions $\lambda(Z_1),\lambda(Z_2)$ are constant on $M$.
Thus, the foliation $\mathcal{F}_{\omega_0|_M}$ coincides with the foliation of the $T^2$-action $\rho$ and therefore, all leaves of $\mathcal{F}_{\omega_0|_M}$ are homeomorphic to $T^2$.
}

Similar to the proof of Theorem \ref{T.main}, we have the following result for arbitrary $r\geq 1$.

\begin{thm}\label{T.main2}
Let $M$ be a $(2n+r)$-dimensional oriented closed $C^{\infty}$-manifold with an
exact presymplectic form $\omega$, that is, there exists a $1$-form $\alpha \in \Omega^1(M)$ such that the $2$-form $\omega=d\alpha$ has rank $2n$ everywhere.
Assume that $M$ admits a locally free $T^r$-action with the following conditions:
\vspace{-9pt}
\begin{enumerate}[{\rm (i)}]
\setlength{\itemsep}{-5pt}
\item The $T^r$-action preserves $\omega$,
\item The functions $\alpha(Z_i)$, $i=2,\ldots,r$, are constant on $M$,
\end{enumerate}\vspace{-9pt}
where $Z_1,Z_2,\ldots,Z_r$ denote the infinitesimal generators of the $T^r$-action.
Then the $r$-dimensional $C^{\infty}$-foliation $\mathcal{F}_{\omega}$ has at least two leaves which are homeomorphic to an $r$-dimensional torus.
Moreover, if the function $\alpha(Z_1)$ is also constant, then $\mathcal{F}_{\omega}$ coincides with the foliation of the $T^r$-action and hence, all leaves of $\mathcal{F}_{\omega}$ are homeomorphic to an $r$-dimensional torus.
\end{thm}

Theorem \ref{T.main2} also gives a partial answer to Question in Section $1$.
If $r=1$, then Theorem \ref{T.main2} agrees with Theorem \ref{T.B-R}.

\rem{
We emphasize that in Theorem \ref{T.main2} (and $2$) we do not assume that $(M,\omega)$ is {\it $r$-contact}, which means that $M$ carries $r$ linearly independent non-vanishing $1$-forms $\alpha_1,\ldots,\alpha_r$ with a splitting
$
TM=\mathcal{R}\oplus
(\cap_i^r \mbox{ker}\, \alpha_i)
$
satisfying that
$
d\alpha_i|_{\cap \mbox{ker}\, \alpha_i}
$
is non-degenerate and
$
\mbox{ker}\, d\alpha_i=\mathcal{R}
$
for every $i$.
In our context, $\mathcal{R}$ corresponds to $\mbox{ker}\, \omega$.
Finamore \cite[Theorem $3.23$]{Finamore} proved that if a closed presymplectic manifold $(M,\omega)$ is $r$-contact with a special metric, then the $r$-dimensional foliation defined by $\mathcal{R}$ has at least two leaves which are homeomorphic to an $r$-dimensional torus.
In the case where $M$ is $r$-contact, by definition, $M$ admits a locally free $\mathbb{R}^r$-action with the infinitesimal generators
$R_1,\ldots,R_r \in \mathcal{R}$.
On the other hand, in Theorem \ref{T.main2}, $M$ does not always satisfy that $Z_1\in \mbox{ker}\, \omega$, as (\ref{eq.Z_1}) in Example \ref{E.1} shows.
Thus, \cite[Theorem $3.23$]{Finamore} and Theorem \ref{T.main2} are independent results.
}

\section*{Acknowledgement}
This paper is a part of the author's Master thesis presented in February 2022 at Ibaraki University.
I would like to thank my supervisor, Prof.\ Hiroshi Iriyeh for suggesting the problem, enlightening discussions and helpful comments.

\Addresses

\end{document}